\documentclass[12pt]{article}
\usepackage{ulem}
\usepackage{graphicx,amsmath,amsfonts,amssymb,color,mathrsfs, amsmath,amsthm}
\usepackage{epsfig}
\newcommand{\chuhao}{\fontsize{19pt}{\baselineskip}\selectfont}
\newcommand{\xiaohao}{\fontsize{10pt}{\baselineskip}\selectfont}

\newcommand{\BOX}{\hfill $\Box$}

 \newcommand{\vc}[1]{{\boldsymbol #1}}

\numberwithin{equation}{section}

 \newtheorem{theorem}{Theorem}[section]
 \newtheorem{lemma}{Lemma}[section]
 
 \newtheorem{proposition}{Proposition}[section]
 
 \newtheorem{remark}{Remark}[section]

 \newtheorem{example}{Example}[section]
 
 \newtheorem{corollary}{Corollary}[section]

\title{\bf\color{black} \chuhao{Perturbation bounds for the stationary distributions of Markov  chains}}
\author{
 Yuanyuan Liu
 \thanks{Corresponding author. Postal address: School
of Mathematics, Railway Campus, Central South University, Changsha,
Hunan, 410075, China; Email address:
liuyy@csu.edu.cn; Tel.: (+86) 731-82655267.}\\
\xiaohao{School of Mathematics, Central South University, Changsha, 410075, China} \\
}
\date{2nd Revision,  July-12, 2012}
\date{July-12, 2012, \ \  Accepted by SIMAX journal}
 
\begin{document}
 \maketitle

\noindent \line(1,0){450}
\\
   \smallskip
\noindent\textbf{Abstract}

In this paper, we are interested in investigating the perturbation bounds
for the stationary distributions for discrete-time or continuous-time Markov chains on a countable state space.
For discrete-time Markov chains, two new norm-wise bounds are obtained. The first bound is rather easy to be obtained  since the needed condition, equivalent to uniform ergodicity, is imposed on the transition matrix directly. The second bound, which holds for a general (possibly periodic) Markov chain, involves finding a drift function. This drift function is closely related with the mean first hitting times. Some $V$-norm-wise bounds are also derived based on the results in \cite{Ka1986}. Moreover, we show how the bounds developed in this paper and one bound given in \cite{Seneta1988} can be extended to continuous-time Markov chains. Several examples are shown to illustrate our results or to compare our bounds with the known ones in the literature.
\medskip

\noindent \textit{AMS Classification:}  60J10; 60J27; 15B51
\medskip

\noindent \textit{Keywords:}  Markov chains; Uniform ergodicity; Stationary distribution; Perturbation theory;
  Mean first hitting times\\
     \smallskip
\noindent \line(1,0){450}\\

\section{Introduction}

\ \ \  \ Let $\Phi(n)$  be a (time-homogeneous) discrete-time  Markov chain (DTMC)  with an irreducible and stochastic transition
matrix $P=(P(i,j))$ on a countable (finite or infinite) state space $\mathbb{E}$. Denote by $P^n=(P^n(i,j))$ the $n$-step transition matrix of $\Phi(n)$. The number $d$, defined by $d=\gcd \{n\geq 1: P^n(i,i)>0\}$ for any (then for all) $i\in \mathbb{E}$,  is called the period for $\Phi(n)$, where $\gcd$ stands for the greatest common divisor.  The chain $\Phi(n)$ chain is said to be aperiodic if $d=1$. Obviously, if $P(i,i)>0$ for some  $i\in \mathbb{E}$, then $\Phi(n)$ is aperiodic. Define $\sigma_C=\inf\{n\geq 0: \Phi(n)\in C\}$ to be the first hitting time on a set $C\subseteq \mathbb{E}$ and write  $m_{ij}=E_i[\sigma_j]$. Note that $m_{ii}=0$.  Let $e$ be  a column vector of all ones.  Suppose that $\Phi(n)$ is perturbed to be another DTMC $\tilde{\Phi}(n)$ with the irreducible and
stochastic transition matrix $\tilde{P}$. Let $\Delta=\tilde{P}-P$.
Suppose that $P$ and $\tilde{P}$ are  positive recurrent
with the unique invariant probability measure (row vector) $\pi$ and $\nu$, respectively.
  Let $\Pi$ be a matrix with equal rows $\pi$. Note that when $\mathbb{E}$ is finite, both $P$ and $\tilde{P}$
  are automatically positive recurrent. We are interested in deriving the perturbation bounds for the difference between $\nu$ and $\pi$ in terms of $\Delta$.

The $V$-norm (see \cite{Ka1986, {Ka1985}}) is introduced as follows.
Let $V$ be a  finite function $V$ on $\mathbb{E}$ bounded away from zero, i.e. $ \inf_{i\in \mathbb{E}} V(i)>0$.
For a finite measure $\mu$, let $\mu (V)=\sum_{i\in \mathbb{E}} \mu(i) V(i)$ and define its $V$-norm to be
  $
   \| \mu \|_V=\sum_{i\in \mathbb{E}} |\mu(i)| V(i).
$
 Let $x$ be a vector on  $\mathbb{E}$ and define its  $V$-norm as
$   \|x\|_V=\sup_{i\in \mathbb{E}} \frac{|x(i)|}{V(i)}. $
 The  $V$-norm for any matrix $L=(L_{ij})$ on $\mathbb{E}\times \mathbb{E}$
 is given by
$
   \|L\|_V=\sup_{i\in \mathbb{E}} \frac{1}{V(i)} \sum_{j\in \mathbb{E}} |L_{ij}| V(j).
$
 When $V\equiv 1$, we omit the subscript $V$ in the notations of $\|
\mu \|_V$, $\|x\|_V$, and  $\|L\|_V$. Note that $\|\mu L x\|_V\leq \|\mu\|_V\|L\|_V\|x\|_V$ and $\|A B\|_V\leq \|A\|_V \|B\|_V$ for any pair of matrices $A$ and $B$ on $\mathbb{E}\times \mathbb{E}$.

We now review some known results on perturbation bounds for a DTMC on a finite state space. The perturbation bounds mainly include the component-wise  bounds  for  $|\nu(k)-\pi(k)|$, $k\in \mathbb{E}$ and the (measure) norm-wise  bounds for
 $\|\nu - \pi\|=\sum_{j\in \mathbb{E}} |\nu(j)-\pi(j)| $. The following two formulas
\begin{equation}
 \nu-\pi = \nu \Delta R,
\end{equation}
\begin{equation}
 \nu-\pi = \nu \Delta  A^{\#},
\end{equation}
derived by \cite{Schweitzer1968} and \cite{Meyer1980}, respectively,  are fundamental for perturbation analysis. Here  $R=(I-P+\Pi)^{-1}$ is the fundamental matrix and $A^{\#}$ is  the group inverse  of $A=I-P$.  The group inverse $B^{\#}=(B^{\#}_{ij})$ of a matrix $B$ is the unique square matrix such that
  \[
    B B^{\#} B=B,\ \  B^{\#} B B^{\#}=B^{\#} \ \ \mbox{and}\ \ B^{\#} B =B B^{\#}.
  \]
 From \cite{Meyer1975}, we know that  $A^{\#}=R-\Pi$. However, $A^{\#}$ has more computational advantages than $R-\Pi$. A lot of component-wise bounds (e.g. \cite{{Cho2000},{Haviv1984},{Funderlic1986}, {Meyer1994}, {Kirkland1998}}) and
norm-wise bounds (e.g. \cite{Seneta1988,{Seneta1991}}) have been obtained in terms of $A^{\#}_{ij}$, $m_{ij}$ or ergodicity coefficient.
For a (possibly negative) matrix $B=(B_{ij})$, the ergodicity coefficient of  $B$ is defined by $\Lambda_1(B)= \frac{1}{2}\sup_{i,j\in \mathbb{E}} \sum_{k\in \mathbb{E}}|B_{ik}-B_{jk}|$. If $\Lambda_1(P)<1$,  Seneta \cite{Seneta1988} derived the following norm-wise bound
 \begin{equation}
\| \nu-\pi\| \leq \frac{\|\Delta\|}{1-\Lambda_1(P)}.
\end{equation}
Subsequently,  Seneta \cite{Seneta1991} obtained another norm-wise bound
\begin{equation}
\| \nu-\pi\| \leq \Lambda_1(A^{\#}) \|\Delta\|,
\end{equation}
which holds even when $\Lambda_1(P)=1$. It was proved by \cite{Kirkland2008} that the bound, given by (1.4), is the smallest (best) norm-wise bound.

It was pointed out in page 13 of \cite{Funderlic1986} that computation of $A^{\#}$ is generally
expensive. Computing the group inverse $A^{\#}$ could be challenging for a large-size finite transition matrix. Moreover, the group inverse $A^{\#}$ can not be used directly for an infinite Markov chain,  since $A^{\#}$ itself needs to be well defined. Thus the bounds characterized by $A^{\#}$ for finite chains can not be simply extended to infinite chains. Let $D=(D_{ij})$ be the deviation matrix defined by $D_{ij}=\sum_{n=0}^\infty (P^n-\Pi)$.  The bound given by (1.3) can be extended to infinite Markov chains whenever $\|D\|<\infty$ (see \cite{Rabta2008}), equivalently, $\Phi(n)$ is uniformly ergodic (see Lemma 2.1 in this paper). This bound is very sensitive when $\Lambda_1(P)$ is close to 1. When $\Lambda_1(P)=1$, we may consider the ergodicity coefficient $\Lambda_1(P^m)$ of the skeleton chain $P^m$ for some positive integer $m$ (see \cite{Mitro2005}), which, however, is not easy to be determined for infinite Markov chains. Hence it is interesting to look for some new bounds which can be expressed in a simple way and can be applied to infinitely countable Markov chains.  Motivated by these issues, we are focused on deriving new norm-wise perturbation bounds for a countable Markov  chain in Section 2. Our approach is based on ergodicity theory. The conditions, imposed on the DTMCs,  are closely related with uniform ergodicity ( i.e. $\|P^n-\Pi\|\rightarrow 0$, as $n\rightarrow \infty$). A simple norm-wise perturbation bound is derived in Section 2.1 by using a ``small set'' condition, which is equivalent to uniform ergodicity, and which is  imposed directly on the transition matrix. This  bound is obtained by bounding the ergodicity coefficient $\Lambda_1(P^m)$ in terms of the ``small set'' condition. Please note that this result holds only for aperiodic Markov chains. Hence we further present another perturbation bound in Section 2.2 for a general (possibly periodic) Markov chain. This bound is given by bounding $R-\Pi$ in terms of the drift condition \textbf{D1$(V,C)$}. When the chain is aperiodic, this drift condition is also equivalent to uniform ergodicity. As a byproduct, a method is proposed to calculate the first hitting times $m_{ij}$, which is different from that in \cite{Hunter2005}.

The more general $V$-norm-wise perturbation bounds for $\| \nu-\pi\|_V$ are developed in the seminal work of \cite{Ka1986,{Ka1985}}. This topic has also gained much interest in the past decades, see, e.g. \cite{Altman, {Heidergott2010}}. This condition $\|R\|_V<\infty$ was used by  \cite{Ka1986} to investigate the $V$-norm-wise perturbation bounds. When $V$ is bounded, $\|R\|_V<\infty$ is equivalent to $\|R\|<\infty$. In this case, the perturbation bounds for $\| \nu-\pi\|_V$ do not make much sense, and we are more interested in the norm-wise bounds. And, the $V$-norm-wise perturbation bounds, obtained by \cite{Ka1986}, are not explicit enough for being used directly to derive the norm-wise bounds. That is why we need to consider the norm-wise bounds separately in Section 2. While $V$ is unbounded, the $V$-norm-wise perturbation bounds enable us to measure the perturbation of the moments of the invariant distribution, which causes essential difference from the norm-wise bounds. In Section 3, we derive some $V$-norm-wise perturbation bounds, expressed in terms of the drift condition \textbf{D2$(V,\lambda, C)$}, for DTMCs based on the results in \cite{Ka1986}.

It is recognized that CTMCs (e.g. \cite{Anderson91}) are important for  modeling
 real phenomena in biology, finance, information, and so on. In the context of queueing theory,
 many queueing models are closely related with CTMCs, for example, the M/M/s/N queue itself is a CTMC. Perturbation analysis is not a new topic here, please see \cite{Tweedie1980} for the component-wise bounds, see \cite{Altman, {Heidergott2010}} for the $V$-norm-wise bounds, and see  \cite{Mitro2006} and \cite{Zeifman1994} for the norm-wise bounds for finite and non-homogeneous CTMCs, respectively. Although a lot of bounds have been developed, it is still worthwhile to develop new and applicable perturbation bounds from different aspects.  In Section 4, we will show how the perturbation bounds developed in  Sections 2 and 3 and the bound given by (1.3) can be extended to CTMCs. Some conclusions are listed in Section 5, and some related results from \cite{Ka1986} are stated in Section 6.

 \section{Norm-wise bounds for DTMCs}

We first extend (1.1) to a countable sate space. To achieve this, we need to define the inverse of the operator $I-P+\Pi$. Let $\ell_1=\{\mu: \|\mu\| <\infty\}$ be the Banach space of all the finite measures. The linear operator $I-P+\Pi: \ell_1\rightarrow \ell_1$  is well defined and its domain is $\ell_1$.   Indeed, for any $\mu \in \ell_1$, $\mu  (I-P+\Pi)=0$ implies  $\mu =0$.  Hence $R=(I-P+\Pi)^{-1}$ is also well defined according to the basic theory about the inverse of a linear operator. Note that $\|R\|$ may be finite or infinite.
Since $\pi R = \pi$ and $\nu \Delta = \nu (I-P+\Pi)-\pi$, we have
\begin{equation}
\nu-\pi = \nu \Delta R = \nu \Delta (R-\Pi).
\end{equation}

If the deviation matrix $D$ exists ( i.e. $D<\infty$), then $\sum_{n=0}^\infty (P-\Pi)^n = D + \Pi $ and

\[
  (I-P+\Pi) \left(\sum_{n=0}^\infty (P-\Pi)^n \right) = \left(\sum_{n=0}^\infty (P-\Pi)^n \right)(I-P+\Pi) = I.
\]
Due to the uniqueness of a linear operator,  we have
\begin{equation}
  R=(I-P+\Pi)^{-1} = D+\Pi.
\end{equation}
From (2.1) and (2.2), we have
\begin{equation}
  \nu-\pi = \nu \Delta (D+\Pi) = \nu \Delta D.
\end{equation}

 Note that if the deviation matrix $D$ exists, then the chain  must be  aperiodic. Based on the above arguments, we can conclude that  (2.3) holds only for an aperiodic Markov chain, while (2.1) holds for a periodic or aperiodic Markov chain.  The following proposition, most parts of which are known (see \cite{Ka1985}), relates  the boundness of $\|R\|$ and $\|D\|$ with the uniform ergodicity.
\begin{proposition}
The following conditions are equivalent to each other:
\begin{description}
\item[(i)] $\Phi(n)$ is aperiodic and $\|R\|<\infty$;

\item[(ii)] $\|D\|<\infty$;

\item[(iii)] $\Phi(n)$ is uniformly ergodic.
\end{description}
\end{proposition}
\proof  It follows from both Theorem 1 and Corollary of Theorem 3 in \cite{Ka1985} that (i) holds if and only if (iii) holds. If (ii)
holds, then the chain must be aperiodic, and
\[
\left\|\sum_{n=0}^\infty (P-\Pi)^n\right\|=\|D+\Pi\|\leq \|D\|+ \|\Pi\|<\infty.
\]
From (2.2), we have $ \|R\| = \|\sum_{n=0}^\infty (P-\Pi)^n\|< \infty$, i.e. (i) holds.  If (iii) holds, then from
Theorem 16.02 in \cite{meyn1993}, we know that there exist positive constants $r<1$ and $b<\infty$ such that
   \begin{equation}
     \|P^n-\Pi \| \leq b r^{- n}
  \end{equation}
for any $n\geq 0$, which implies that $\Phi(n)$ is aperiodic and
$
    \|D\|\leq \sum_{n=0}^\infty \|(P^n-\Pi)\|< \infty,
$
 i.e. (ii) holds.  \BOX

\subsection{A norm-wise bound based on uniform ergodicity}

To derive the main results in this subsection, we need the concept of a small set. Let ${\cal B}(\mathbb{E})$ be the set composed of  all the subsets of $\mathbb{E}$. A set $C$ is called a small set if
there exist a positive integer $m$ and a non-trivial measure $\nu_m$ on ${\cal B}(\mathbb{E})$ such that
\[
   P^m(i, B):=\sum_{k \in B} P^m(i,k) \geq \nu_m(B)
\]
for any $i\in C$ and any $B\in {\cal B}(\mathbb{E})$. For DTMCs on a countable state space, every finite set is a small set. It is known from Theorem 16.0.2 in \cite{meyn1993} that $\Phi_n$ is uniformly ergodic if and only if the whole state space $\mathbb{E}$ is a $\nu_m$-small set for some $m$.

\begin{theorem}
If the  state space $\mathbb{E}$ is $\nu_m$-small for some positive integer $m$ and some non-trivial measure $\nu_m$, then
\begin{equation}
 \| \nu-\pi\| \leq \frac{\|P^m - \tilde{P}^m \|}{\nu_m({\mathbb{E}})} \leq  \frac{m}{\nu_m({\mathbb{E}})}\|\Delta\|.
\end{equation}
In particular, if there exists some positive integer $m$ such that $\sum_{k\in \mathbb{E}}\delta_m(k)>0$, where $\delta_m(k):=\inf_{i\in \mathbb{E}}P^m (i, k)$ ,  then   (2.5) holds with $\nu_m(\mathbb{E})=\sum_{k\in \mathbb{E}}\delta_m(k)$.
\end{theorem}

\begin{lemma} If $\Lambda_1(P^m)<1$ for some $m\geq 1$, then
\begin{equation*}
  \| \nu-\pi\| \leq \frac{\|P^m - \tilde{P}^m \|}{1-\Lambda_1(P^m)}.
\end{equation*}
\end{lemma}
\proof  Let $\{\Phi(nm), n\geq 0\}$ be the $m$-skeleton chain of $\Phi(n)$. Then the chain $\Phi(nm)$ has one-step transition probability matrix $P^m$. Since $\pi P =\pi$, we have $\pi P^m=\pi$, which implies that $\pi$ is also the  invariant distribution of the chain $\Phi(nm)$. Similarly, we know that $\nu$ is also the invariant distribution of the $m$-skeleton chain $\tilde{\Phi}(nm)$. The assertion is obtained immediately by applying (1.3) to  both skeleton chains $\Phi(nm)$ and $\tilde{\Phi}(nm)$.

\begin{remark}  Following the arguments in Section 2 of \cite{Seneta1988}, we can extend this lemma easily to more general $p$-norm $\|\cdot\|_p$, $1\leq p\leq \infty$ (the $\ell_p$ norm on the space of real row vectors). Here, we only state this simple case (i.e. $p=1$), which was first presented in \cite{Mitro2005} using different arguments,  to avoid introducing too many mathematical notations.
\end{remark}
\begin{lemma} If  the  state space $\mathbb{E}$ is $\nu_m$-small for some positive integer $m$ and some non-trivial measure $\nu_m$, then
\begin{equation*}
\Lambda_1(P^m)\leq 1-\nu_m(\mathbb{E}).
\end{equation*}
\end{lemma}
\proof Since the  state space $\mathbb{E}$ is $\nu_m$-small, we have
\[
  P^m(i, k) \geq \nu_m(k),\ \ i\in \mathbb{E}
\]
for any fixed $k\in \mathbb{E}$. Hence, for any $i,k$, we can find a non-negative real number $d_{ik}$ such that
\[
  P^m(i,k) = \nu_m(k) + d_{ik},
\]
which implies that
\begin{eqnarray*}
  \Lambda_1(P^m)&=& \frac{1}{2}\sup_{i,j\in \mathbb{E}} \sum_{k\in \mathbb{E}}|P^m(i, k)-P^m(j, k)|\nonumber\\
  &\leq& \frac{1}{2}\sup_{i,j\in \mathbb{E}} \sum_{k\in \mathbb{E}}(d_{ik}+d_{jk})\nonumber\\
  &=& 1-\sum_{k\in \mathbb{E}} \nu_m(k)\nonumber\\
  & =& 1-\nu_m(\mathbb{E}).
\end{eqnarray*}

\noindent\textit{Proof of Theorem 2.1} \ \ By Lemmas 2.1 and  2.2, we have
\begin{equation}
   \| \nu-\pi\| \leq \frac{\|P^m - \tilde{P}^m \|}{\nu_m({\mathbb{E}})},
\end{equation}
which is the first inequality of (2.5). It is easy to derive
\begin{equation}
\|P^m - \tilde{P}^m \|\leq \|(P-\tilde{P})\| \|P^{m-1}+P^{m-2}\tilde{P}+\cdots+ P\tilde{P}^{m-2}+\tilde{P}^{m}\|\leq m \|\Delta\|.
\end{equation}
From (2.6) and (2.7), we obtain the second inequality of (2.5).

To prove the second part of this assertion, define a set function $\nu_m$ as follows
\[
  \nu_m(C)=\sum_{k\in C} \nu_m(k),\ \ C\in {\cal B}(\mathbb{E}),
\]
where $  \nu_m(k)=\nu_m(\{k\}) = \delta_m(k),  k\in \mathbb{E}.$
Obviously, the non-negative set function $\nu_m$ constitutes a non-trivial measure on ${\cal B}(\mathbb{E})$.
Observe that for any $A\in {\cal B}(\mathbb{E})$
\[
 \inf_{i\in \mathbb{E}} P^m (i, A) = \nu_m(A),
\]
from which, and the first assertion, we obtain the second part of the theorem.
\BOX
\begin{remark}
 We could have obtained the following perturbation bound more directly at the cost of a worse bound that is twice as big as the one given by (2.5). Suppose that the  state space $\mathbb{E}$ is $\nu_m$-small. Then Theorem 16.2.4 in \cite{meyn1993}shows that
         \begin{equation*}
         \|P^n -\Pi\|\leq 2(1-\nu_m(\mathbb{E}))^{\lfloor \frac{n}{m}\rfloor},
         \end{equation*}
where $\lfloor \frac{n}{m} \rfloor$ denotes the greatest integer not exceeding $\frac{n}{m}$.
From (2.3) and Proposition 2.1, we have
\[
\| \nu-\pi\| \leq \|D\| \|\Delta\| \leq   \sum_{n=0}^\infty \|P^n-\Pi \|  \|\Delta\| \leq  \frac{2 m}{\nu_m(\mathbb{E})}\|\Delta\|.
\]
\end{remark}

 To apply this result, it is helpful to know which kind of Markov chains could be uniformly ergodic. When the state space $\mathbb{E}$ is finite, an irreducible, aperiodic and  positive recurrent DTMC $\Phi(n)$ is always uniformly ergodic. For a Markov chain on an infinite state space, uniform ergodicity usually requires that the state space should have a ``central state'', which is accessible from all other states in finite time. To see this, we note that a positive recurrent Markov chain is uniformly ergodic if and only if for any fixed $j\in \mathbb{E}$, there is an integer $N$ such that $\inf_{i\in \mathbb{E}} P^n(i, j)\geq \frac{\pi_j}{2}>0$ for all $n\geq N$. Also, we know from Section 2 in \cite{Hou2004} that  $\Phi(n)$ cannot be uniformly ergodic if  $P$ is a Feller transition matrix, i.e. $\lim_{i\rightarrow \infty} P(i,j)=0$ for any fixed $j\in \mathbb{E}$. These observations give us some insight into uniform ergodicity for DTMCs on an infinite state space.

\begin{example} Consider the DTMC on $\mathbb{E}=\mathbb{Z}_+$ with the following lower-Hessenberg transition matrix:
\begin{equation*}
 P=\left (
 \begin{array}{lllll}
   b_0 & a_0 & 0& 0 & ... \\
   b_1 & a_1 & a_0 & 0& ... \\
   b_2& a_2 & a_1 & a_0 & ... \\
   b_3& a_3& a_2 & a_1 & ... \\
   ... &... & ... & ...& ... \
   \end{array}\right ), \end{equation*}
where $\mathbb{Z}_+$ is the set of all non-negative integers. Suppose that the chain is
irreducible and $\sum_{k=0}^\infty a_k < 1$. It is known that the transition matrix of the embedded GI/M/1 queue with negative arrivals is of the above structure, in which  $a_j, j\geq 0$ take specific forms such that $\sum_{k=0}^\infty a_k < 1$. Since $\sum_{k=0}^\infty a_k < 1$, we have $d_0(1)=\inf_{i\in \mathbb{E}} P(i,0)\geq 1-\sum_{k=0}^\infty a_k$. From Theorem 2.1 we have the norm-wise perturbed bound
 $ \|\Delta\| \frac{1}{1-\sum_{k=0}^\infty a_k}. $
 \BOX
  \end{example}

\begin{example}  Consider the DTMC on $\mathbb{E}=\mathbb{Z_+}$ with the following  transition matrix elements:
\begin{equation*}
 P(i,j)=\left \{
  \begin{array}{lllll}
   q, \ \ \mbox{if}\ \ j=i+1, \ \ i\geq 0,\\
   p, \ \ \mbox{if} \ \ j=0, \ \ i= 0,\\
 p, \ \ \mbox{if}\ \ j=0, \ \ \mbox{i is odd},\\
 p, \ \ \mbox{if}\ \ j=1, \ \ i\geq 1 \ \ \mbox{and i is even},\\
  0, \ \ \mbox{else},
    \end{array}\right.
  \end{equation*}
   where $p$ and $q$ are positive numbers such that $p+q=1$.
Calculating the elements in the  first column of $P^2$, we
have $d_0(2)=\inf_{i\geq 0} P^2 (i, 0)= p^2$. Thus we obtain the perturbation bound
$\frac{2}{p^2} \| \Delta\|$ from Theorem 2.1. \BOX
\end{example}

 To compare our bound with the best one (1.4), we borrow two examples from the literature. The first one is from \cite{Funderlic1986}, which  models the mammillary systems in compartmental analysis, and which was used to compare the perturbation bounds in \cite{Kirkland2008}.
\[
 P=\left (
 \begin{array}{ccccccccccc}
 0.74 & 0.11 & 0 & 0 & 0 & 0 & 0 & 0.15\\
0 & 0.689 & 0 & 0 & 0.011 & 0&  0&  0.3\\
0&  0 & 0 & 0.4 & 0 & 0 & 0 & 0.6\\
0 & 0 & 0 & 0.669 & 0.011 & 0 & 0 & 0.32\\
0 & 0&  0 & 0 & 0.912 & 0 & 0 & 0.088\\
0 & 0 & 0 & 0 & 0 & 0.74 & 0 & 0.26\\
0 & 0 & 0 & 0 & 0 & 0 & 0.87 & 0.13\\
0.15 & 0 & 0.047 & 0 & 0 & 0.055 & 0.27 & 0.478
       \end{array}\right )
  \]
According to \cite{Kirkland2008}, we know $\Lambda_1(A^{\#})=11.3352$. From this and (1.4), we get the perturbation bound $11.3352 \|\Delta\|$. 
Observing the last column of $P$ and using Theorem 2.1, we have
the slightly bigger bound $ \frac{1}{0.088} \|\Delta\|=11.3636\|\Delta\|$.

The second one is from \cite{Meyer1975}, whose transition matrix is given by
\[
P= \frac{1}{4}\left (
 \begin{array}{cccccc}
  0 & 2 & 2 & 0 \\
 2 & 0 & 2&  0 \\
2 &  1 & 0 & 1 \\
1 & 1 & 1 & 1
       \end{array}\right ).
\]
According to \cite{Meyer1975}, $A^{\#}$ is given as follows, from which and (1.4), we obtain the perturbation bound  $1.5512 \|\Delta\|$. Computing $P^2$ (given below), we obtain from Theorem 2.1 the bound $3.2 \|\Delta\|$.
\[
       A^{\#}= \frac{2}{1083}\left (
 \begin{array}{cccccc}
  265 & -61 & -96 & -108 \\
 -96 & 300 & -96 &  -108 \\
 -115 & -137 & 246 & 6 \\
 -210 & -156 & -210 & 576
       \end{array}\right ),\ \
         P^2= \frac{1}{16}\left (
 \begin{array}{cccccc}
  8 & 2 & 4 & 2 \\
 4 & 6 & 4&  2 \\
3 &  5 & 7 & 1 \\
5 & 4 & 5 & 2
       \end{array}\right ).
\]

\subsection{A norm-wise bound based on a drift condition}

We have known that a uniformly ergodic Markov chain is necessarily aperiodic. The bound given by Theorem 2.1 holds only for aperiodic chains. In this subsection, we will make use of the formula (2.1) and the following drift condition to derive a norm-wise perturbation bound for a general Markov chain which is possibly periodic.

 \textbf{D1($V,C$)}: There exist a bounded non-negative
function $V$ and a finite set $C$ such that
\begin{equation}\label{}
 \left \{
 \begin{array}{cc}
  \sum\limits_{j\in \mathbb{E}}P(i,j)V(j)\leq V(i)-1,\ \  & i\notin C,\\
V(i)=0, \ \ \ & i\in C.
 \end{array}\right .
\end{equation}

 It is well known (e.g. \cite{liu2006}) that  an aperiodic and irreducible Markov chain is uniformly ergodic if and only if its transition matrix $P$ satisfies the drift condition. As will be shown in the following proposition, this drift function is greater than or equal to a sequence of the first hitting times. A more general form of the following proposition is proposed by Theorem 2.1 in \cite{liu2011}.

 \begin{proposition}  Let $C$ be any fixed finite set in $\mathbb{E}$.
The function $V$, defined by $V(i)=E_i [\sigma_C], i\in \mathbb{E}$,
is the minimal non-negative solution to (2.8), and  $V$ satisfies (2.8) with equality. Note
 that the minimal solution means that if there is another solution $\tilde{V}$, then we always have
 $V(i)\leq \tilde{V}(i)$ for any $i\in \mathbb{E}$.

\end{proposition}

  \begin{theorem}  Let $i_0$ be any fixed state in $\mathbb{E}$.
If  $P$ satisfies \textbf{D1($V,C$)} for $C=\{i_0\}$,  then
  \begin{equation}
   \| \nu-\pi\| \leq 2 \left(\sup_{i\in \mathbb{E}} V(i)\right)^2 \|\Delta\|,
   \end{equation}
and
   \begin{equation}
   \| \nu-\pi\| \leq 2 \inf_{i_0\in \mathbb{E}}\left(\sup_{i\in \mathbb{E}} m_{ii_0}\right)^2 \|\Delta\|.
   \end{equation}
  \end{theorem}

  \proof  (i) Let $g$ be an indicator function on $\mathbb{E}$ given by
  \[
g(i)=\left \{
  \begin{array}{lllll}
1, \ \ & \mbox{if}\ \  i=i_0,\\
0 ,\ \ & \mbox{else}\ \ i\neq i_0.\\
\end{array}\right.
\]
 Define a measure $\alpha$ on ${\cal B}(\mathbb{E})$ by
 \[
 \alpha(j)= P(i_0, j), j\in \mathbb{E} \ \ \mbox{and}\ \ \alpha(A)=\sum_{i\in A} \alpha(i), A\in {\cal B}(\mathbb{E}).
 \]
 Let $T=(T_{ij})$ be the matrix  given by
\begin{equation*}
T_{ij}=P(i,j)-g(i) \alpha(j)=\left \{
  \begin{array}{lllll}
0, \ \ & \mbox{if}\ \  i=i_0,\\
 P(i, j) ,\ \ & \mbox{if}\ \ i\neq i_0.
\end{array}\right.
\end{equation*}
 For any $\varepsilon>0$, define the sequence $\{\hat{V}(i), i\in \mathbb{E}\}$ by $\hat{V}(i)=V(i)+\varepsilon$, $i\in \mathbb{E}$. Since $V$ satisfies (2.8),  we have
  \begin{equation*}
     \sum_{j\in \mathbb{E}} T_{ij}  \hat{V}(j) \leq \left(1-\frac{1}{\sup_{i\in \mathbb{E}}{\hat{V}(i)}}\right) \hat{V}(i)
  \end{equation*}
  for any $i\in \mathbb{E}$, which follows
  \begin{equation}
    \|T\|_{\hat{V}}\leq 1-\frac{1}{\sup_{i\in \mathbb{E}}{\hat{V}(i)}}.
  \end{equation}
  It is easy to derive
 \[
   \|P\|_{\hat{V}}\leq \max\left\{\frac{\sup_{i\in \mathbb{E}}\hat{V}(i)}{\varepsilon}, 1-\frac{1}{\sup_{i\in \mathbb{E}}{\hat{V}(i)}}\right\}<\infty.
  \]
  Hence the condition in Theorem 2 in \cite{Ka1986} (see condition (ii) in Proposition 6.1 in the Appendix) is satisfied. It follows from (6) in \cite{Ka1986} (see (6.1) in the Appendix), we have
  \begin{equation*}
    R-\Pi=\Pi \left(\pi \sum_{n=0}^\infty T^n  e I -\sum_{n=0}^\infty T^n \right) + \sum_{n=0}^\infty T^n(I-\Pi),
  \end{equation*}

  Using (2.1) and the fact that $\Delta \Pi=0$, we have
  \begin{equation}
\nu-\pi =  \nu \Delta (R-\Pi) = \nu \Delta \sum_{n=0}^\infty T^n(I-\Pi).
\end{equation}
   From (2.11),  we obtain
  \begin{eqnarray*}
    \left\| \sum_{n=0}^\infty T^n(I-\Pi) \right\|&\leq&  \| I-\Pi\|  \left\|\sum_{n=0}^\infty T ^n \right\|  \nonumber\\
    &\leq&  2 \sup_{i\in \mathbb{E}}\hat{V}(i) \left\|\sum_{n=0}^\infty T ^n \right\|_{\hat{V}}\nonumber\\
    &\leq& 2 \sup_{i\in \mathbb{E}}\hat{V}(i) \sum_{n=0}^\infty (\|T\|_{\hat{V}})^n \nonumber\\
    &\leq& 2 \left(\sup_{i\in \mathbb{E}}\hat{V}(i)\right)^2.
  \end{eqnarray*}
   Since $\varepsilon$ can be given arbitrarily, we have
   \begin{equation}
     \left\| \sum_{n=0}^\infty T^n(I-\Pi) \right\| \leq 2 \lim_{\varepsilon \downarrow 0}\left(\sup_{i\in \mathbb{E}} V(i)+\varepsilon\right)^2 \leq 2 \left(\sup_{i\in \mathbb{E}} V(i)\right)^2.
   \end{equation}
 The perturbation bound (2.9) follows from (2.12) and (2.13) immediately.

 (ii) From Proposition 2.1, we know that $\sup_{i\in \mathbb{E}} m_{ii_0}\leq \sup_{i\in \mathbb{E}}{V(i)}<\infty$, and that the sequence of $\{m_{ii_0}, i\in \mathbb{E}\}$ also satisfies the drift condition.  It is well known that for an irreducible Markov chain, $\sup_{i\in \mathbb{E}} m_{ii_0}<\infty$  for some state $i_0\in \mathbb{E}$ if and only if $\sup_{i\in \mathbb{E}} m_{ii_0}<\infty$ for any sate $i_0\in \mathbb{E}$. By (2.9), we have
 \begin{equation}
   \| \nu-\pi\| \leq 2 \left(\sup_{i\in \mathbb{E}} m_{ii_0}\right)^2 \|\Delta\|,
   \end{equation}
   which follows the bound (2.10) since the state $i_0$ is taken arbitrarily.
  \BOX

 To apply this result, a key point is to find a drift function, which is the usual way to verify ergodicity in the context of ergodic theory. For finite Markov chains, we can solve (2.8) with equality for $C=\{i_0\}$, i.e. to compute the first hitting times $m_{ii_0}$. This is always feasible since it is equivalent to solving a finite system of linear equations. This way of computing $m_{ij}$ is different from the one proposed in \cite{Hunter2005}, where $m_{ij}$ is calculated in terms of a $g$-inverse $G$ of $I-P$ (i.e. $(I-P)G(I-P)=I-P$). For infinite Markov chains, we do not have such a general procedure. However, we can make some suggestions on this. One way is still to solve  (2.8) with equality by making use of the structure of the transition matrix and the ``minimal nonnegative'' property, which will be illustrated by the following Example 2.4.  The other way is to construct a drift function such that (2.8) holds with inequality instead of strict equality. As will be shown by the last two examples at the end of this section, it is possible to construct a simple constant solution for special models.

 \begin{example}
Consider the DTMC  on  $\mathbb{E}:=\{0, 1, \cdots, n\}$ with the following birth-death matrix:
  \[
      P =   \left (
      \begin{array}{lllllll}
      c_0& b_0 & 0& \cdots & 0& 0&  0 \\
      a_1 & c_1& b_1& \cdots& 0& 0&  0\\
      0 & a_2  & c_2& \cdots& 0 &0&  0 \\

      \vdots & \vdots & \vdots  & \ddots& \vdots &\vdots& \vdots \\
       0 & 0  &0  &\cdots& c_{n-2}&b_{n-2}&  0 \\
       0 &  0 & 0  & \cdots & a_{n-1} &c_{n-1} & b_{n-1} \\
       0 &  0 & 0  & \cdots &0& a_{n} & c_n \\
      \end{array}  \right),
  \]
where the coefficients $a_i, b_i$ and $c_i$ are such that $P$ is stochastic and irreducible. Note that $P$ is periodic with periodicity $d=2$ if $c_i=0$ for all $0\leq i\leq n$. For any fixed $j\in \mathbb{E}$, we define the drift function $V$ by  $V(i)= m_{ij}, i\in \mathbb{E}$.
Substituting the value of
$P(i, j)$ into (2.8) with equality and inducing on $m$ gives
\begin{eqnarray*}
  V(m)-V(m+1) &=&\frac{a_n}{b_n}[V(m-1)-V(m)]+\frac{1}{b_m}
 =\cdots
 =\frac{1}{b_m \mu(m)}  \sum_{k=0}^m  \mu(k)
\end{eqnarray*}
 for any $m$, $0\leq m\leq {j-1}$, where $V(-1)=0$,   and the sequence $\{\mu(k), 0\leq k\leq n\}$  is given by
 \[
   \mu(0)=1,\ \ \mu(k)=\frac{b_0\cdots b_{n-1}}{a_1\cdots, a_n}, 1\leq k\leq n.
 \]
Summing over $m$ from $i$ to $j-1$ yields
\begin{equation}
 m_{ij}= V(i)=\sum_{m=i}^{j-1} \frac{1}{b_m\mu(m)} \sum_{k=0}^m \mu(k), \ \
  i<j.
\end{equation}
Then we consider the case of $i>j$. Substituting the value of
$P(i,j)$, we have
\[
V(n)-V(n-1) =\frac{1}{a_n},
\]
and
\[
 V(m+1)-V(m) =\frac{a_m}{b_m}(V(m)-V(m-1))-\frac{1}{b_m},\ \  j+1\leq  m \leq n-1,
 \]
 where $V(j)=0$.
 Solving these equations gives
\begin{equation}
 m_{ij}= V(i)= \sum_{m=j}^{i-1} \frac{1}{b_m \mu(m)} \sum_{\ell=m+1}^n \mu(\ell),\ \ j< i\leq n.
\end{equation}
From (2.10), (2.15) and (2.16), we obtain the perturbation bound
\[
 \| \nu-\pi\| \leq 2 \inf_{0\leq i_0\leq n} \left(\max \left\{\sum_{m=i_0}^{n-1} \frac{1}{b_m \mu(m)} \sum_{\ell=m+1}^n \mu(\ell), \sum_{m=0}^{i_0-1} \frac{1}{b_m \mu(m)} \sum_{\ell=0}^m \mu(\ell) \right\}\right)^2  \|\Delta\|,
\]
where we make the convention that $\sum_{k=n}^{n-1} a_k = \sum_{k=0}^{-1} a_k =0$ for any sequence of $a_k$.
\BOX
 \end{example}

\begin{example}
Consider the DTMC on $\mathbb{E}=\mathbb{Z}_+$ with  the following  transition matrix elements:
  \[
    P(i,j)=\left\{
   \begin{array}{lllll}
 1, &  i=0, j=1 \\
  p_i, &  i\geq 1, j=0 \\
   q_i, &  i\geq 1, j=i+1 \\
    0, &  \mbox{else}. \\
       \end{array}\right.
  \]
where $p_i$ and $q_i$ are positive numbers such that $p_i+q_i=1$.  Define the drift function $V$ by  $V(i)= m_{i0}, i\in \mathbb{E}$.  Taking $i_0=0$ in (2.8) with equality, we have
 \begin{equation}
   V(n+1)= \frac{V(n)}{q_n} - \frac{1}{q_n} =\cdots = \frac{V(1)}{\prod_{k=1}^n q_k} - \sum_{j=0}^{n-1} \frac{1}{\prod_{k=0}^j q_k},\ \ n\geq 1.
 \end{equation}
Since $V$ is the minimal and non-negative solution, we  have
\begin{equation}
V(1)=\sup_{n\geq 1}  \prod_{k=1}^n q_k \sum_{j=0}^{n-1} \frac{1}{\prod_{k=0}^j q_{n-k} } =  \sup_{n\geq 1} \sum_{j=0}^{n-1} \frac{1}{\prod_{k=j+1}^n q_{n-k} }.
\end{equation}
Hence $m_{i0}, i\geq 1$ are obtained by (2.17) and (2.18). In particular, for $p_i=p$ and $q_i=q=1-p$,  we have
$ m_{i 0}=\frac{1}{1-q} -\frac{1}{q^i}$ for any $i\geq 1$. Using (2.9), we obtain the perturbation bound $\frac{2}{p^2}$.
\BOX

\end{example}

 To show how to find a drift function such that (2.8) holds with inequality, we consider Example 2.1 and Example 2.2. For both examples, we choose $i_0=0$.

For Example 2.1, we let $V(0)=0$ and $V(i)=\frac{1}{1-\sum_{k=0}^\infty a_k}$. Then  for $i\geq 1$
\[
 \sum_{j=0}^\infty P(i,j) V(j) = \frac{\sum_{k=0}^{i} a_k}{1-\sum_{k=0}^\infty a_k}=\frac{1}{1-\sum_{k=0}^\infty a_k}- \frac{1-\sum_{k=0}^{i} a_k}{1-\sum_{k=0}^\infty a_k} \leq V(i)-1.
\]
By (2.9), we obtain the perturbation bound  $ \frac{2}{(1-\sum_{k=0}^\infty a_k)^2} \|\Delta\|.$ This bound is worse than the one given by Theorem 2.1.

 The transition matrix $P$ in Example 2.2 is modified into a periodic one with periodicity $d=2$, by changing two elements $P(0,0)=p$ and $P(0,1)=q$ into $0$ and $1$ respectively, and keeping all the other elements unchanged.  For any fixed positive number $\lambda$, let $V(0)=0$, $V(1)=\frac{\lambda}{p}$ and $V(i)=\frac{1+\lambda}{p}, i\geq 2$. Then we have
\[
  \sum_{j=0}^\infty P(i,j) V(j) = \left\{
   \begin{array}{ccc}
   \frac{1+\lambda}{p} q = \frac{1+q\lambda}{p}-1=V(1)-1 , &  i=1, \\
 \frac{1+\lambda}{p} q <\frac{1+\lambda}{p}-1=V(i)-1, &  \ \  i\geq 2,  i \  \mbox{is odd}, \\
 \frac{\lambda}{p}p + \frac{1+\lambda}{p} q = \frac{1+\lambda}{p}-1=V(i)-1, &  i\geq 2, i \  \mbox{is even}. \\
       \end{array}\right.
\]
Letting $\lambda\downarrow 0$, we obtain the perturbation bound $\frac{2}{p^2}\|\Delta\|$ from Theorem 2.2.

\section{$V$-norm-wise bounds  for DTMCs}

To consider the $V$-norm-wise bounds  for DTMCs, we need to change (2.1) into
\begin{equation}
  \nu = \pi (I- \Delta (R-\Pi))^{-1}=\pi \sum_{n=0}^\infty [\Delta (R-\Pi)]^n,
\end{equation}
under the assumption that  $\|\Delta (R-\Pi)\|_V<1$. The $V$-norm-wise bounds were first investigated by \cite{Ka1986}  through a detailed analysis of the boundness $\|R\|_V$ and an explicit expression of $R-\Pi$.  Based on Corollary 2 in \cite{Ka1986} (see Proposition 6.2 in the Appendix), we obtain the perturbation bounds for $ \| \nu-\pi\|_V$ in terms of the following drift condition.

\textbf{D2($V,\lambda, b,C$)}: There exists a finite function $V$ bounded away from zero, some finite set $C$,
and positive constants $\lambda<1, b<\infty$ such
that
\begin{equation}\label{1.3}
\sum_{j\in \mathbb{E}} P(i,j)V(j) \leq \lambda V(i) + b I_C(i),\ \ \ \
i\in \mathbb{E}.
\end{equation}

\begin{corollary}
Let $i_0$ be any fixed state in $\mathbb{E}$. Suppose that $P$ satisfies \textbf{D2($V,\lambda, b, C$)} for $C=\{i_0\}$.
\begin{description}
\item[(i)] Let  $c=1+\|e\|_V \|\pi\|_V$. If $\|\Delta\|_V<\frac{1-\lambda}{c}$, then $\tilde{\Phi}(n)$ is positive recurrent and
\begin{equation}
  \| \nu-\pi\|_V \leq   \frac{c \|\pi\|_V \|\Delta\|_V}{1-\lambda-c\|\Delta\|_V}.
\end{equation}

\item[(ii)] If $V\geq 1$ and $\|\Delta\|_V<\frac{(1-\lambda)^2}{b+1-\lambda}$,  then $\tilde{\Phi}(n)$ is positive recurrent and
\begin{equation}
  \| \nu-\pi\|_V \leq    \frac{b(b+1-\lambda)\|\Delta\|_V}{(1-\lambda)^3 -(1-\lambda)(b+1-\lambda)\|\Delta\|_V}.
\end{equation}
\end{description}
\end{corollary}
\proof (i) Let  $T=(T_{ij})$ be exactly  the same matrix as that defined in the proof of Theorem 2.2, i.e. $T$ is formed from $P$ by changing $P(i_0, j), j\in \mathbb{E}$ into $0$s and keeping the other elements unchanged.
Obviously, all the three conditions of Corollary 2 in \cite{Ka1986} are satisfied. So we have the first assertion immediately from  Corollary 2 in \cite{Ka1986}.

 (ii) Multiplying both side of (3.2) by $\pi(i)$ and summing over $i$,  we obtain
\[
   \sum_{i\in \mathbb{E}} \pi(i) \sum_{j\in \mathbb{E}} P(i,j)V(j) \leq \lambda \pi(V) + b \pi(i_0).
\]
Using the invariance of $\pi$ ( i.e. $\pi P =\pi$) derives
\[
   \|\pi\|_V =\pi(V)\leq \frac{ b \pi(i_0)}{1-\lambda}\leq \frac{b}{1-\lambda}.
\]
We obtain the second assertion immediately  from the first one, by replacing $\|\pi\|_V$ and $c$ in the condition of (i) and in  (3.3) with their upper bounds $\frac{b}{1-\lambda}$ and $1+\frac{b}{1-\lambda}$, respectively. \BOX

\begin{remark} (i) It seems impossible to define a  measure $\alpha$ and a function $h$ so that the single point set $\{i_0\}$ in this corollary  can be changed into a  finite set $C$. (ii) The second bound is completely dependent on the parameters $\lambda, b, V$ and $\|\Delta\|_V$ at the cost of decreasing the accuracy of the bound. To investigate a specific model, we should try to use the first assertion whenever $\pi(V)$ can be computed.
\end{remark}

\begin{remark}
This drift condition \textbf{D2($V,\lambda, b,C$)} is sufficient and necessary for an aperiodic chain $\Phi_n$
to be $V$-uniformly ergodic, i.e. $\|P^n-\Pi\|_V\rightarrow 0$ as $n\rightarrow \infty$. When $V$ is bounded, the condition \textbf{D2($V,\lambda, b,C$)} is theoretically equivalent to the condition \textbf{D1($V,C$)}. However, it is much harder to decide the former, since it involves two more parameters. That is the reason why   we choose the latter to derive the norm-wise bound in Section 2. An aperiodic chain $\Phi_n$ is  geometrically ergodic but not uniformly ergodic if and only if $P$ satisfies \textbf{D2($V,\lambda, b,C$)} for an unbounded function $V$. The $V$-norm-wise perturbation bounds really make sense when $V$ is unbounded.
\end{remark}

 \section{Perturbation bounds for CTMCs}

  We now show how the discrete-time results developed thus far may be lifted through the $h$-approximation chain to obtain
analogous results for CTMCs.

Let $\Phi_t$ be a CTMC on a countable state
space $\mathbb{E}$ with an
irreducible and conservative intensity matrix $Q=(Q_{ij})$. Throughout this section, we assume that $Q$ is uniformly bounded (i.e. $\sup_{i\in \mathbb{E}} Q_i<\infty$, where $Q_i:=-Q_{ii}$). Let $P^t=(P^t(i,j))$ be the corresponding unique $Q$-function. It is known that $P^t = e^{Q t}=\sum_{n=0}^\infty \frac{(Q t)^n}{n!}$.  We suppose
that $\Phi_t$ is positive recurrent with the unique invariant probability measure $\pi$. The
matrix $Q$  is perturbed to be another irreducible and uniformly bounded
 intensity matrix $\tilde{Q}=(\tilde{q}_{ij})$. The
corresponding $\tilde{Q}$-process $\tilde{\Phi}_t$ is assumed to be positive recurrent
with the invariant probability measure $\nu$. Let $\Delta=\tilde{Q}-Q$.  We are interested in investigating the perturbation bounds for $\nu-\pi$ in terms of $\Delta$.

We now introduce the $h$-approximation chain (e.g. \cite{Anderson91}), which is a crucial technique adopted to extend the perturbation bounds from  DTMCs to CTMCs. For the uniformly bounded intensity matrix $Q$, let $h< {(\sup_{i\in \mathbb{E}} Q_i)}^{-1}$ be the length of the time discretisation
interval. The transition probabilities $P^h(i,j)$ for $\Phi(t)$ have first
order approximations $P_h(i,j)=(I+hQ)_{ij}, \ \ i,j\in \mathbb{E}$. The  matrix $P_h$ is stochastic, irreducible and aperiodic. The chain $\Phi_h(n)$, with $P_h=(P_h(i,j))$ as its one-step transition matrix, is called  the $h$-approximation
chain of $\Phi(t)$.  Define $D=(D_{ij})=\int_0^\infty
(P^t-\Pi)dt$ to be the deviation matrix of $\Phi_t$ and let $D_h$ be the deviation matrix of $P_h$.  We have the following basic relations between $\Phi(t)$ and $\Phi_h(n)$:
\begin{description}
\item[(i)] $\Phi(t)$ and $\Phi_h(n)$  have the same invariant probability measures;

\item[(ii)] $\Phi(t)$ and $\Phi_h(n)$ are equivalent in positive recurrence, geometric ergodicity, and uniform ergodicity ( see, e.g. \cite{liu2006, {Spiek91}});

\item[(iii)] $D=\frac{1}{h}D_h$ (see, e.g. \cite{Van Doorn2003}).
\end{description}

\subsection{Norm-wise bounds for CTMCs}

To consider the perturbation of the invariant probability measure for a CTMC, we can equivalently consider the perturbation of the invariant probability measure for its $h$-approximation chain.  For any  $h<\left[\max\{\sup_{i\in \mathbb{E}}Q_i, \sup_{i\in \mathbb{E}} \tilde{Q}_i\}\right]^{-1},$ let $P_h$ and $\tilde{P}_h$ be the $h$-approximation chains of $\Phi_t$ and $\tilde{\Phi}_t$, respectively.
Since $P_h$ is aperiodic, $D_h$ can be used to characterize the perturbation bounds for the $h$-approximation chain, which enables us to investigate the perturbation bounds for CTMCs only in terms of  the deviation matrix $D$.

 The following drift condition, which is sufficient and necessary for $\Phi(t)$ to be uniformly ergodic (i.e. $\|P^t-\Pi\|\rightarrow 0$), can be found in Chapter 6 of \cite{Anderson91}.

\textbf{D1'($V,C$)}:  There exist a bounded non-negative
function $V$ and a finite set $C$ such that
\begin{equation}\label{}
 \left \{
 \begin{array}{cc}
  \sum\limits_{j\in \mathbb{E}}Q_{ij}V(j)\leq -1,\ \  & i\notin C,\\
V(i)=0, \ \ \ & i\in C.
 \end{array}\right .
\end{equation}

We are now in  a position to state the norm-wise perturbation bounds for CTMCs.

\begin{theorem}  Let $i_0$ be any fixed sate in $\mathbb{E}$.
\begin{description}
\item[(i)] If the $Q$-process $\Phi_t$ is uniformly ergodic, then we have
\begin{equation}
\| \nu-\pi\| \leq \| D\| \|\Delta\|.
\end{equation}
\item[(ii)]  If $\Lambda_1(Q)>0$, then
  \begin{equation}
 \| \nu-\pi\| \leq \frac{\Delta}{\Lambda_1(Q)},
 \end{equation}
where $\Lambda_1(Q)=  \frac{1}{2}\inf_{i,j}
\left[|Q_{ii}-Q_{ji}|+|Q_{ij}-Q_{jj}|-\sum_{s\neq i, s\neq j}
|Q_{is}-Q_{js}|\right]$.
   \item[(iii)] Let $\delta_k=\inf_{i\neq k} Q_{ik}$. If $\sum_{k\in \mathbb{E}} \delta_k>0$, then $\Phi_t$ is uniformly ergodic, and
   \begin{equation}
    \| \nu-\pi\| \leq \frac{1}{\sum_{k\in \mathbb{E}} \delta_k} \|\Delta\|.
   \end{equation}
 \item[(IV)]  If Q satisfies \textbf{D1'($V,C$)} for $C=\{i_0\}$, then we have
 \begin{equation}
  \| \nu-\pi\| \leq  2 \left( \sup_{i\in \mathbb{E}}
  V(i) \right)^2 \|\Delta\|.
  \end{equation}
 \end{description}
\end{theorem}
\proof  (i) Since $\Phi(t)$ is uniformly ergodic,  $\Phi_h(n)$ is also uniformly ergodic. It implies from Proposition 2.1 that $\|D_h\|<\infty$. From (2.3), we have
 \begin{equation}
   \nu-\pi\ = \nu \cdot \Delta_h\cdot  D_h = \nu \Delta D,
 \end{equation}
where  $\Delta_h=\tilde{P}_h-P_h$. We then have (4.2) from (4.6).

(ii) Since $\Lambda_1(Q)<0$, we can choose small enough $h$ such that
\begin{eqnarray}
  \Lambda_1(P_h)&=& \frac{1}{2}\sup_{i,j\in \mathbb{E}} \left[|1+h Q_{ii}-h Q_{ji}|+ |h Q_{ij}-1-hQ_{jj}|+\sum_{s\neq i, s\neq j}
h |Q_{is}-Q_{js}| \right]\nonumber \\
  &=& \frac{1}{2}\sup_{i,j\in \mathbb{E}} \left[2-h|Q_{ji}-Q_{ii}|-h| Q_{ij}-Q_{jj}|+ \sum_{s\neq i, s\neq j}
h|Q_{is}-Q_{js}| \right]\nonumber \\
  &=& 1-h\Lambda_1(Q)<1.
  \end{eqnarray}
Hence we have (4.3) from (1.3) and (4.7).

(iii) Choose small enough $h$ such that $h\leq \frac{1}{ \sup_{i\in \mathbb{E}} Q_i+ \sum_{k\in \mathbb{E}}\delta_k}$. Then we have
\[
  \inf_{i\in \mathbb{E}} P_h(i,k)\geq h \delta_k
\]
for any $i\in \mathbb{E}$,  from which and (2.5), we obtain
\[
   \| \nu-\pi\| \leq  \frac{1}{h \sum_{k\in \mathbb{E}}\delta_k}\|\Delta(h)\| = \frac{1}{ \sum_{k\in \mathbb{E}}\delta_k}\|\Delta\|.
\]

(IV) Since  Q satisfies \textbf{D1'($V,C$)} for $C=\{i_0\}$, we can obtain
\begin{equation}
  \sum_{j\in \mathbb{E}} Q_{ij} \hat{V}(j) \leq - \hat{\lambda} \hat{V}(i), \ \ i\neq i_0,
\end{equation}
where $\hat{\lambda}=\frac{1}{\sup_{i\in \mathbb{E}} \hat{V}(i)}$ and $\hat{V}(i)=V(i)+\varepsilon$ for any $i\in \mathbb{E}$, where $\varepsilon$ is an arbitrarily given positive number.
Transferring (4.8) to the $h$-approximation chain gives
\[
  \sum_{j\in \mathbb{E}} P_h(i,j) \hat{V}(j) \leq (1-\hat{\lambda} h) \hat{V}(i), \ \  i\neq i_0.
\]
Similar to the proof of Theorem 2.2, we can obtain the bound (4.5). \BOX

\begin{remark}  If the state space $\mathbb{E}$ is finite, then $Q$ is  uniformly bounded and  the  $Q$-process is  uniformly  ergodic. Hence all the conditions (i)-(iii) hold automatically .
\end{remark}

\subsection{$V$-norm-wise bounds for CTMCs}

We will use the following drift condition, which is equivalent to $V$-uniform ergodicity for a CTMC, to find the $V$-norm-wise bounds for CTMCs. These bounds parallel to the ones in Corollary 3.1.

 \textbf{D2'($V,\lambda, b,C$)}: There exists a finite function $V$ bounded away from zero, some finite set $C$,
and positive constants $\lambda, b<\infty$ such
that
\begin{equation}\label{1.3}
\sum_{j\in \mathbb{E}} Q_{ij}V(j)\leq -\lambda V(i) + b I_C(i),\ \ \ \
i\in\mathbb{E}.
\end{equation}

\begin{theorem} Let $i_0$ be any fixed state in $\mathbb{E}$. Suppose that $Q$ satisfies \textbf{D2'($V,\lambda, b, C$)} for $C=\{i_0\}$.
\begin{description}
\item[(i)] Let  $c=1+\|\vc{e}\|_V \|\pi\|_V$. If $\|\Delta\|_V<\frac{\lambda}{c}$, then $\tilde{\Phi}(t)$ is positive recurrent, and
\begin{equation}
  \| \nu-\pi\|_V \leq   \frac{c \|\pi\|_V \|\Delta\|_V }{\lambda-c\|\Delta\|_V}.
\end{equation}

\item[(ii)] If $V\geq 1$ and $\|\Delta\|_V<\frac{\lambda^2}{b+\lambda}$, then
\begin{equation}
  \| \nu-\pi\|_V \leq    \frac{b(b+\lambda)\|\Delta\|_V}{\lambda^3 - \lambda(b+\lambda)\|\Delta\|_V}.
\end{equation}
\end{description}
\end{theorem}
\proof Since $Q$ is bounded and $Q$ satisfies \textbf{D2'($V,\lambda, b, \{i_0\}$)}, we have
\[
  \sum_{j\in \mathbb{E}}P_h(i,j) V (j) \leq (1-\lambda h) V(i) + bh I_{\{i_0\}}(i).
\]
The bounds (4.10) and (4.11) follow from (3.3) and (3.4), respectively.
\BOX

\begin{remark} Suppose that $Q$ satisfies \textbf{D2'($V,\lambda, b, C$)} for $V\geq 1$ and $C=\{i_0\}$.  From (3.1) and the above arguments, we know that if  $\|\Delta\|_V<\frac{\lambda}{c}$ or $\|\Delta\|_V<\frac{\lambda^2}{b+\lambda}$, then
$
   \nu = \pi  \sum_{n=0}^\infty [\Delta D]^n.
$
In particular, let $\tilde{Q}=Q+ \varepsilon G$. Assume that $G e=0$ and  $\|G\|_V\leq g_1<\infty$. Then for any $\varepsilon$ such that $\varepsilon<\frac{\lambda}{c g_1}$ or $\varepsilon<\frac{\lambda^2}{(b+\lambda) g_1}$, we have
$
   \nu = \pi  \sum_{n=0}^\infty \varepsilon ^n (G D)^n.
$
This  result is better than Theorem 2.1 in \cite{Altman}, since the convergence domain  of $\varepsilon$ for this series expansion is more computable and accurate.
\end{remark}

\begin{example}
Consider the  CTMC on $\mathbb{E}=\mathbb{Z}_+$ with the following intensity matrix:
\begin{equation*}
 Q=\left (
 \begin{array}{ccccc}
   a_0 & a_1 & a_2 & a_3 & ... \\
   b_0 & b_1 & b_2 & b_3 & ... \\
   0   &  b_0 &b_1 & b_2 & ... \\
   0   &  0  & b_0   & b_1& ... \\
   \vdots  &\vdots & \vdots & \vdots & \vdots \
   \end{array}\right ),
  \end{equation*}
where $\{a_i, i\in \mathbb{Z}_+\} and \{b_i,i\in \mathbb{Z}_+\}$ are two sequences of real
numbers such that $Q$ is stable and conservative. A Markovian queue with batch arrivals is a particular case of this chain. For a sequence of real numbers $\{c_k, k\in \mathbb{Z_+}\}$, define $C(z)=\sum_{i=0}^\infty c_k z^k$ to be the generating function, and let $\phi_C$ be the radius of convergence of $C(z)$. It is known (see, e.g. \cite{liu2006}) that $\Phi_t$  is ergodic if and only if $ \sum_{k=0}^\infty k a_k<\infty$ and  $B'(1-) <0$ (i.e. $\sum_{k=1}^\infty kb_{k+1}<b_0$).
 Suppose that this chain is ergodic, $\phi=\min\{\phi_A, \phi_B\}>1$ and $0<B(\phi_B)<\infty$.
Define $\rho=\sup\{z: B(z) \leq  0\}$. Since $B''(z)>0$, $B(z)$ is a convex function in $[0, \rho]$. Thus the function $\frac{-B(z)}{z}$ is continuous
in $[1, \rho]$. Hence $\hat{\lambda}:=\max\left\{\frac{-B(z)}{z}:z\in [1,\rho]\right\}$ can be attained at some point $z_0\in [1, \rho]$. Furthermore, we have $z_0>1$ and $\hat{\lambda}>0$, since $B(1)=0$ and  $B'(1-)<0$.  Let $V(i)=z_0^i, i\geq 0$. We have
\[
  \sum_{j\in \mathbb{E}}Q_{ij} V (j) = -\frac{-B(z_0)}{z_0} V(i)=-\hat{\lambda} V(i), \ \ i\geq 1,
\]
and $   \sum_{j\in \mathbb{E}}Q_{0j} V (j) = A(z_0)<\infty$. Thus the equality in (4.9) holds for $i_0=0$, $\lambda=\hat{\lambda}$ and $b=A(z_0)+\hat{\lambda}$, which implies $\pi(V)=b/\hat{\lambda}$. By (4.10) or (4.11), we have
\begin{equation}
   \| \nu-\pi\|_V \leq   \frac{b(b+\hat{\lambda})\|\Delta\|_V}{\hat{\lambda}^3 - \hat{\lambda}(b+\hat{\lambda})\|\Delta\|_V}.
\end{equation}

In particular, let $a_1=\sigma$, $a_i=0, i\geq 2$; $b_0=\mu$, $b_2=\sigma$ and $b_i=0, i\geq 3$. Then the chain $\Phi(t)$ becomes the well-known M/M/1 queue. If $\Phi(t)$ is positive recurrent, i.e. $\sigma<\mu$, then we have $z_0=\sqrt{\mu/\sigma}$, $V(i)=z_0^i,$ $\hat{\lambda}=(\sqrt{\sigma}-\sqrt{\mu})^2$ and $b=\mu- \sqrt{\mu \sigma}$. Hence the $V$-norm-wise bound can be given explicitly by (4.12).
 \BOX
 \end{example}

\section{Concluding remarks}

Two new  perturbation bounds for DTMCs are derived by giving computable values of $\ell$ in $\| \nu-\pi \|\leq \ell \|\Delta\|$. The $V$-norm-wise perturbation bounds are also considered. These bounds developed for DTMCs are further extended to CTMCs. If a norm-wise bound is such that $\ell \|\Delta\|\geq 2$, then the bound is entirely useless, since we automatically have $\|\nu-\pi\|\leq 2$. Hence $\|\Delta\|$ is usually assumed to be small.

 A DTMC $\Phi(n)$  is said to be strongly $V$-stable (strongly stable for $V\equiv 1$), if every stochastic transition matrix $\tilde{P}$ in $\{\tilde{P}:\|\Delta\|_V:=\|P-\tilde{P}\|_V<\varepsilon\}$ has a unique invariant probability measure $\nu$ such that $ \|\nu-\pi\|_V\rightarrow 0$ as $\|\Delta\|_V\rightarrow 0$. It is known from Theorem 1 in \cite{Ka1986} (see (i) of Proposition 6.1 in the Appendix) that $\Phi(n)$ is $V$-stable if and only if $\|R\|_V<\infty$.  It can be easily seen that $\Phi(n)$ is strongly stable under the condition in Theorem 2.1 or Theorem 2.2, and that $\Phi(n)$ is $V$-strongly stable under the condition in Corollary 3.1.

Similarly, we can define strong stability for CTMCs.  A chain $\Phi(t)$ is is said to be strongly $V$-stable (strongly stable for $V\equiv 1$), if every conservative intensity matrix $\tilde{Q}$ in $\{\tilde{Q}:\|\Delta\|_V:=\|Q-\tilde{Q}\|_V<\varepsilon\}$ has a unique invariant probability measure $\nu$ such that $ \|\nu-\pi\|_V\rightarrow 0$ as $\|\Delta\|_V\rightarrow 0$.  Obviously, $\Phi(t)$ is  strongly stable under the condition in Theorem 4.1, and $\Phi(t)$ is $V$-strongly stable under the condition in Theorem 4.2.

Currently, we only consider perturbation bounds for CTMCs with the uniformly bounded intensity matrices. When the intensity matrices are unbounded,  the $h$-approximation chain method can not be used, and the bounds given in Section 4 may fail to hold. It is meaningful to know which bounds still hold. To investigate this issue requires some different methods, which is a topic for future research.
\section{Appendix}

The following proposition, taken from Theorems 1 and 2 in \cite{Ka1986}, are important for investigating the strong stability and perturbation bounds in $V$-norm.

 \begin{proposition} Let $V$ be a finite function bounded away from zero.
The DTMC $\Phi(n)$ is strongly $V$-stable  if and only if either of the following condition holds:
\begin{description}
\item[(i)]
 $\|P\|_V<\infty$ and $\|R\|_V=\|(I-P+\Pi)^{-1}\|_V<\infty$.

\item[(ii)] (1) $\|P\|_V<\infty$; (2) There are some finite non-negative measure $\alpha$
and some non-negative bounded function $g$ on $\mathbb{E}$ such that $T=(T_{ij})$ is a non-negative matrix, where $T_{ij}=P(i, j)-g(i) \alpha(j)$;
and (3) there exist positive numbers $\lambda<1$ and $m\geq 1$ such that $T^m V(x)\leq \lambda V(x)$.

\end{description}
\end{proposition}

The following Proposition is taken from formula (6) and  Corollary 2 in \cite{Ka1986}.

\begin{proposition}
Suppose that condition (ii) in Proposition 6.1 holds. Then we have
\begin{equation}
R-\Pi=\Pi \left[\pi (I-T)^{-1} e I -(I-T)^{-1} \right] + (I-T)^{-1}(I-\Pi),
\end{equation}
where $(I-T)^{-1}=\sum_{n=0}^\infty T^n$. Furthermore, if $\|\Delta\|_V<\frac{1-\lambda}{c}$ with $c=1+\|e\|_V \|\pi\|_V$, then
\[
  \| \nu-\pi\|_V \leq   \frac{c \|\Delta\|_V \|\pi\|_V}{1-\lambda-c\|\Delta\|_V}.
\]

\end{proposition}

\begin{center}
\bigskip \noindent{\bf\large Acknowledgements}
\end{center}

The author is deeply grateful to two anonymous referees for  very constructive comments and suggestions, which have allowed the author to
improve the presentation of this paper significantly. The author acknowledges the support from  D\'{e}partement d'Informatique at the
Universit\'{e} Libre de Bruxelles during the visit during which this research was completed. This work was
supported in part by the Fundamental Research Funds for the Central Universities (grant number 2010QYZD001) and the National Natural Science Foundation of China (grant numbers  10901164, 11071258).

  \end{document}